\documentclass[12pt]{article} %fleqn
 \usepackage{amsmath}
 \usepackage{amssymb}
 \usepackage{dsfont}
 \usepackage{amsthm}
   \font\twelvebm                       = cmmib10 at 12truept
   \font\tenbm                          = cmmib10 at 10truept
   \font\sevenbm                        = cmmib10 at 7truept
\tenbm \scriptscriptfont9              = \sevenbm

  at 10truept
 at 10truept  at 10truept  at 10truept 
at 10truept
\mathchardef \BGamma            = "0900 \mathchardef
\BDelta = "0901 \mathchardef \BTheta            = "0902
\mathchardef \BLambda           = "0903 \mathchardef \BXi
= "0904 \mathchardef \BPi               = "0905
\mathchardef \BSigma = "0906 \mathchardef \BUpsilon
= "0907 \mathchardef \BPhi = "0908 \mathchardef \BPsi
= "0909 \mathchardef \BOmega = "090A \mathchardef \Balpha
= "090B \mathchardef \Bbeta             = "090C
\mathchardef \Bgamma = "090D \mathchardef \Bdelta
= "090E \mathchardef \Bepsilon = "090F \mathchardef \Bzeta
= "0910 \mathchardef \Beta = "0911 \mathchardef \Btheta =
"0912 \mathchardef \Biota = "0913 \mathchardef \Bkappa
= "0914 \mathchardef \Blambda = "0915 \mathchardef \Bmu
= "0916 \mathchardef \Bnu = "0917 \mathchardef \Bxi
= "0918 \mathchardef \Bpi = "0919 \mathchardef \Brho
= "091A \mathchardef \Bsigma            = "091B
\mathchardef \Btau = "091C \mathchardef \Bupsilon
= "091D \mathchardef \Bphi = "091E \mathchardef \Bchi
= "091F \mathchardef \Bpsi = "0920 \mathchardef \Bomega
= "0921 \mathchardef \Bvarepsilon       = "0922
\mathchardef \Bvartheta         = "0923 \mathchardef
\Bvarpi = "0924 \mathchardef \Bvarrho = "0925 \mathchardef
\Bvarsigma = "0926 \mathchardef \Bvarphi           = "0927
\mathchardef \bA        = "0941 \mathchardef \bB        =
"0942 \mathchardef \bC        = "0943 \mathchardef \bD
= "0944 \mathchardef \bE        = "0945 \mathchardef \bF
= "0946 \mathchardef \bG        = "0947 \mathchardef \bH
= "0948 \mathchardef \bI        = "0949 \mathchardef \bJ
= "094A \mathchardef \bK        = "094B \mathchardef \bL
= "094C \mathchardef \bM        = "094D \mathchardef \bN
= "094E \mathchardef \bO        = "094F \mathchardef \bP
= "0950 \mathchardef \bQ        = "0951 \mathchardef \bR
= "0952 \mathchardef \bS        = "0953 \mathchardef \bT
= "0954 \mathchardef \bU        = "0955 \mathchardef \bV
= "0956 \mathchardef \bW        = "0957 \mathchardef \bX
= "0958 \mathchardef \bY        = "0959 \mathchardef \bZ
= "095A \mathchardef \ba        = "0961 \mathchardef \bb
= "0962 \mathchardef \bc        = "0963 \mathchardef \bd
= "0964
\mathchardef \bee       = "0965 %%%I CHANGED IT FROM \be; SEE IN THE BEGGINING.
\mathchardef \bff       = "0966 \mathchardef \bg        =
"0967 \mathchardef \bh        = "0968
\mathchardef \bj        = "096A \mathchardef \bk        =
"096B \mathchardef \bl        = "096C \mathchardef \bm
= "096D \mathchardef \bn        = "096E \mathchardef \bo
= "096F \mathchardef \bp        = "0970 \mathchardef \bq
= "0971 \mathchardef \br        = "0972 \mathchardef \bs
= "0973 \mathchardef \bt        = "0974 \mathchardef \bu
= "0975 \mathchardef \bv        = "0976 \mathchardef \bw
= "0977 \mathchardef \bx        = "0978 \mathchardef \by
= "0979 \mathchardef \bz        = "097A

\font\tencb            = cmssbx10 scaled \magstep4
\font\eigcb = cmssbx10 scaled \magstep2 \textfont8
= \tencb \scriptfont8           = \eigcb \scriptscriptfont8
= \eigcb \mathchardef\bAs       = "1841
\def\Asem#1#2{\mathop{\vrule height10.5pt depth5.5pt width0pt\bAs}_{#1}^{#2}}
\def\asem#1#2{
          \ifmmode
         \ifinner
            \raise0.9pt\hbox{$\scriptstyle\bAs$}_{#1}^{#2}
         \else
            \Asem{#1}{#2}
         \fi
          \fi
          }

\newtheorem{theo}{\small\bf Theorem}[section]

\newtheorem{lem}{\small\bf Lemma}[section]

\newtheorem{prop}{\small\bf Proposition} %[section]

\newtheorem{rem}{\small\bf Remark}[section]

\newtheorem{defi}{\small\bf Definition}[section]

\newtheorem{cor}{\small\bf Corollary}[section]

\newtheorem{example}{\small\bf Example}[section]

\newcommand{\be}{\begin{equation}}
\newcommand{\ee}{\end{equation}}
\newcommand{\E}{\mathds{E}}
\newcommand{\Var}{\mbox{\rm \hspace*{.2ex}Var\hspace*{.2ex}}}
\newcommand{\Cov}{\mbox{\rm \hspace*{.2ex}Cov\hspace*{.2ex}}}

\newcommand{\R}{\mathds{R}}

 \newcommand{\lead}{\mbox{\rm{lead}}}
 \newcommand{\xwrosA}{$\begin{array}{c}\vspace{-2ex}\\}
 \newcommand{\xwrosB}{\\ \vspace{-2ex}\end{array}$}

 \newcommand{\supp}{\mbox{supp}}

 \title{\bf Some counterexamples concerning maximal correlation and linear
 regression\footnote{Work partially supported by the University of
 Athens Research Grant 70/4/5637.}}

 \author{\large
 Nickos Papadatos\footnote{\mbox{
 e-mail:\ {\tt npapadat@math.uoa.gr}
 \  \
 url:\ {\tt http://users.uoa.gr/$\sim$npapadat/}}}
 }

 \date{\small
 Section of Statistics and O.R.,
 Department of Mathematics,
 University of Athens,
 \\ Panepistemiopolis, 157 84 Athens, Greece
 \vspace{-2em}
 }

 \newcommand{\lt}{}

\begin{document}

 \maketitle

 %\vspace*{-2em}

 \thispagestyle{empty}

 \begin{abstract}
 \noindent
 A class of examples concerning the
 relationship of linear
 regression and maximal
 correlation is provided.
 More precisely, these examples show that
 if two random variables have (strictly) linear
 regression on each other, then their
 maximal correlation is not necessarily equal to
 their (absolute) correlation.
 %becomes maximal.
 \end{abstract}
 {\footnotesize
 %
 %{\it MSC}:  Primary 62G30; Secondary 62E10, 60E15.
 %
 %\noindent
 }
 %Characterization of Exponential Distribution; Order Statistics; Records; Splitting model.}

 %\vspace{-1em}

 \section{Maximal correlation and linear regression}
 \label{sec.1}
 Let $(X,Y)$ be a bivariate random vector such that its
 Pearson correlation coefficient,
 \be
 \label{eq.rho}
 \rho(X,Y):=\frac{\Cov(X,Y)}{\sqrt{\Var(X)}\sqrt{\Var(Y)}},
 \ee
 is well defined. If $W$ is a non-degenerate random variable
 then $L_2^*(W)$ is defined to be the class of measurable
 functions $g:\R\to\R$ such that
 $0<\Var[g(W)]<\infty$. Under the present notation, the maximal
 correlation coefficient is defined as (Gebelein, 1941; Hirschfeld, 1935)
 \be
 \label{eq.maxcorr}
 R(X,Y):=\sup_{g_1\in L_2^*(X),
 \ g_2\in L_2^*(Y)}
 \rho(g_1(X),g_2(Y)).
 \ee
 Due to results of Sarmanov (1958a, 1958b), it was believed
 for some time that if
 both $X$ and $Y$ have
 linear regression on each other, i.e., if
 for some constants $a_0$, $a_1$, $b_0$,
 $b_1$,
 \be
 \label{eq.linear}
 \E(X|Y)=a_1 Y+a_0
 \ \mbox{(a.s.)},
 \ \ \
 \E(Y|X)=b_1 X+b_0
  \ \mbox{(a.s.)},
 \
 \ee
 then
 \be
 \label{eq.R=r}
 R(X,Y)=|\rho(X,Y)|.
 \ee
 %see Rohatgi and Sz\'{e}kely (1992).
 The implication (\ref{eq.linear})$\Rightarrow$(\ref{eq.R=r})
 was cited in a number of subsequent works
 related to maximal correlation of order
 statistics and records,
 including Rohatgi and Sz\'{e}kely (1992),
 Arnold, Balakrishnan and Nagaraja (1998, p.\ 101),
 Sz\'{e}kely and Gupta (1998), David and
 Nagaraja (2003, p.\ 74),
 Ahsanullah (2004, p.\ 23)
 and Barakat (2012).
 However, as we shall show below, this implication
 is not valid
 even in the case of a strictly linear regression,
 $a_1 b_1\neq 0$.
 %It seems
 %that
 %%%%this implication
 %it
 %is based on
 %an incorrect argument given in the proof of
 %Sarmanov's Theorem 2 (1958a, 1958b), saying that
 %``if there is a monotone function in the spectrum
 %of a stochastic (correlation) kernel then it always
 %belongs to the first (largest) eigenvalue''; this error
 %may goes back to Sarmanov (1946) or Nomokonov (1950).
 Note that if $R(X,Y)>0$ then the converse implication,
 (\ref{eq.R=r})$\Rightarrow$(\ref{eq.linear}),
 is valid; see R\'{e}nyi (1959, p.\ 447)
 and Dembo, Kagan and Shepp (2001).

 Examples of uncorrelated random
 variables $X$, $Y$ with
 (trivial) linear regression
 \be
 \label{eq.trivregr}
 \E(X|Y)=\E(Y|X)=0  \ \ \mbox{(a.s.)}
 \ee
 and $R(X,Y)>0=|\rho(X,Y)|$ are known for a long time.
 For instance, P.\ B\'{a}rtfai has calculated
 $R(X,Y)=1/3$ for a uniform in the interior of
 the unit disc. This result was extended by
 P.\ Cs\'{a}ki and J.\ Fischer for the
 %case where
 %$(X,Y)$ is
 uniform distribution in the domain
 $|x|^p+|y|^p<1$ ($p>0$), in which case
 $R(X,Y)=(p+1)^{-1}$;
 see R\'{e}nyi (1959, p.\ 447) and
 Cs\'{a}ki and Fischer (1963). Furthermore,
 Sz\'{e}kely and M\'{o}ri (1985) extended
 this result to the multivariate case and with different
 exponents. Moreover, in response to a question asked by
 Sid Browne of Columbia University, Dembo, Kagan and Shepp (2001)
 constructed a pair
 %of uncorrelated random variables
 $(X,Y)$ satisfying
 %$\E(X|Y)=\E(Y|X)=0$
 (\ref{eq.trivregr})
 and $R(X,Y)=1$.
 (Observe that the same is true
 %when $(X,Y)$ is
 for the uniform distribution
 in the four-point domain
 $\{(0,\pm1),(\pm1,0)\}$.)
 % and  $\rho(X,Y)=0$,
 %thus showing that (\ref{eq.linear}) does not
 %imply (\ref{eq.R=r}).
 Using characterizations of
 Vershik (1964) and Eaton (1986),
 they also showed that for
 any non-Gaussian spherically symmetric random
 vector $(U_1,\ldots,U_k)$, with covariance
 matrix of rank $\geq 2$, there exists
 a pair of
 uncorrelated linear
 forms,
 \[
 X=a_1U_1+\cdots+a_k U_k,
 \ \ Y=b_1U_1+\cdots+b_k U_k,
 \]
 such that (\ref{eq.trivregr}) is fulfilled
 %$X$ and $Y$ have linear regression on each other
 and $R(X,Y)>|\rho(X,Y)|=0$.

 However, in the author's opinion,
 it is important to definitely know
 that (\ref{eq.linear}) does not imply
 (\ref{eq.R=r}) even in the non-trivial
 linear regression case.
 Indeed, if this implication were valid
 in the particular case where
 $a_1b_1\neq 0$,
 then several works concerning
 %Moreover, the
 characterizations of distributions
 through maximal correlation
 of order statistics and records --
 including the papers by Terrell (1983),
 Sz\'{e}kely and M\'{o}ri (1985), Nevzorov (1992),
 L\'{o}pez-Bl\'{a}zquez and Casta\~{n}o-Mart\'{i}nez
 (2006), Casta\~{n}o-Mart\'{i}nez, L\'{o}pez-Bl\'{a}zquez
 and Salamanca-Mi\~{n}o (2007), Papadatos and Xifara (2012)
 -- would be reduced to trivial consequences
 of this implication. The same is
 true for the main result in Dembo, Kagan
 and Shepp (2001), since it is easily checked that
 for the partial sums $S_k=X_1+\cdots+X_k$, based on an iid
 sequence with mean $\mu$ and finite non-zero variance,
 \[
 \E(S_{n+m}|S_n)=S_n+ m\mu \ \ \mbox{(a.s.),} \ \ \
 \E(S_{n}|S_{n+m})=\frac{n}{n+m}S_{n+m} \ \ \mbox{(a.s.)}.
 \]

 The purpose of the present note is to present a
 quite general
 class of random vectors $(X,Y)$, with $X$ and $Y$
 possessing strictly linear
 regression on each other, and such that
 $R(X,Y)>|\rho(X,Y)|>0$.
 This class is elementary and
 it is defined in the next section.

 \section{Counterexamples}
 Let $f_1$ and $f_2$ be
 two univariate probability densities (with respect
 to Lebesgue measure on $\R$)
 with bounded supports,
 $\supp(f_i)\subseteq[\alpha_i,\omega_i]$,
 %and $\supp(f_2)\subseteq[\alpha_2,\omega_2]$,
 $-\infty<\alpha_i<\omega_i<\infty$ ($i=1,2$).
 It is well known that there exists a uniquely defined
 orthonormal polynomial system $\{\phi_n(x)\}_{n=0}^{\infty}$,
 standardized by $\lead(\phi_n):=p_n>0$, where
 $\lead(\phi_n)$ denotes the principal coefficient of
 $\phi_n$. Also, there exists a uniquely defined
 orthonormal polynomial system $\{\psi_n(x)\}_{n=0}^{\infty}$,
 standardized by $\lead(\psi_n):=q_n>0$. Each system is
 complete in the corresponding $L_2$-space, since
 for any real $t$,
 \[
 \int_{-\infty}^{\infty} e^{tx}f_1(x) dx<\infty
 \ \
 \mbox{and}
 \ \
 \int_{-\infty}^{\infty} e^{ty}f_2(y) dy<\infty;
 \]
 see, e.g., Berg and Christensen (1981) or
 Afendras, Papadatos and Papathanasiou (2011).
 Since every polynomial is uniformly bounded in any
 finite interval, we can find
 constants $c_n$, $d_n$ such that
 \[
 1<\sup_{\alpha_1\leq x \leq \omega_1}|\phi_n(x)|=c_n<\infty,
 \ \
 1<\sup_{\alpha_2\leq y\leq \omega_2}|\psi_n(y)|=d_n<\infty,
 \ \ \
 n=1,2,\ldots \ .
 \]
 Consider an arbitrary real
 sequence $\{\rho_n\}_{n=1}^{\infty}$ such
 that
 \be
 \label{eq.rho_n}
 \sum_{n=1}^{\infty}{|\rho_n|} c_nd_n\leq 1,
 \ee
 e.g., $\rho_n=6(\pi^2n^2 c_n d_n)^{-1}$ ($n=1,2,\ldots$)
 or $\rho_n=\lambda n$ ($n=1,\ldots,N$)
 and $\rho_n=0$, otherwise, where $0<\lambda\leq (\sum_{n=1}^N nc_nd_n)^{-1}$.
 Then, the function
 \be
 \label{eq.density}
 f(x,y):=f_1(x)f_2(y)\left(1+\sum_{n=1}^{\infty}
 \rho_n \phi_n(x)\psi_n(y) \right),
 \ \
 (x,y)\in [\alpha_1,\omega_1]\times[\alpha_2,\omega_2],
 \ee
 and $f:=0$ outside
 $[\alpha_1,\omega_1]\times[\alpha_2,\omega_2]$,
 is a bivariate probability density with marginal densities
 $f_1$, $f_2$; this is so because, due to (\ref{eq.rho_n}),
 the series in (\ref{eq.density}) converges, for each
 $(x,y)$ in the domain of definition,
 to a value greater than or equal to $-1$.
 (Actually, the series converges uniformly  and
 absolutely in
 $[\alpha_1,\omega_1]\times[\alpha_2,\omega_2]$.)
 \
 Therefore, $f(x,y)$ is nonnegative.
 %, and is defined to be zero
 %outside $[\alpha_1,\omega_1]\times[\alpha_2,\omega_2]$.
 Next, it is easily checked that its integral over
 $\R^2$ equals 1,
 due to the orthonormality of the polynomials. Finally, it
 is obvious that the marginal densities of $f$ are $f_1$,
 $f_2$.

 Assume now that the random vector $(X,Y)$ has density $f$.
 Then $X$ has density $f_1$ and $Y$ has density $f_2$. Moreover,
 versions of the conditional densities are given by
 \begin{eqnarray*}
 f_{X|Y}(x|y)
 &\hspace{-1ex}=\hspace{-1ex}&
 f_1(x)\left(1+\sum_{n=1}^{\infty} \rho_n
 \phi_n(x)\psi_n(y)\right),
 \ \
 \alpha_1\leq x \leq \omega_1
 \ \
 (\mbox{for each } y\in\supp(f_2)),
 \\
 f_{Y|X}(y|x)
 &\hspace{-1ex}=\hspace{-1ex}&
 f_2(y)\left(1+\sum_{n=1}^{\infty} \rho_n
 \phi_n(x)\psi_n(y)\right),
 \ \
 \alpha_2\leq y\leq \omega_2
 \ \
 (\mbox{for each } x\in\supp(f_1)).
 \end{eqnarray*}
 Due to the orthonormality of the polynomials
 it follows that for all $n\geq 1$,
 \be
 \label{eq.linregr}
 \E(\phi_n(X)|Y)=\rho_n \psi_n(Y)
 \ \ \mbox{(a.s.)},
 \ \ \
 \E(\psi_n(Y)|X)=\rho_n \phi_n(X)
 \ \ \mbox{(a.s.)}.
 \ee
 Clearly, if $\rho_1\neq 0$,
 (\ref{eq.linregr}) shows that $X$ and $Y$ have
 strictly linear regression on each other. In fact, it is
 easily checked, using (\ref{eq.linregr}) and induction on $n$,
 that
 \be
 \label{eq.polreg}
 \E(X^n|Y)=\frac{\rho_n q_n}{p_n} Y^n +P_{n-1}(Y)
 \ \ \mbox{(a.s.)},
 \ \ \
 \E(Y^n|X)=\frac{\rho_n p_n}{q_n} X^n +Q_{n-1}(X)
 \ \ \mbox{(a.s.)},
 \ee
 where $P_{n-1}(t)$ and $Q_{n-1}(t)$ are polynomials
 of degree at most $n-1$ in $t$.
 %(the polynomials may be different
 % in the above conditional expectations).
 Using (\ref{eq.polreg}) and the main result of
 Papadatos and Xifara (2012),
 or directly from (\ref{eq.linregr}),
 it is a simple matter to conclude that
 $R(X,Y)=\sup_{n\geq 1} |\rho_n|$.
 Since the choice of $\{\rho_n\}_{n=1}^{\infty}$
 is quite arbitrary
 (see (\ref{eq.rho_n})),
 it follows that
 \[
 R(X,Y)>|\rho(X,Y)|=|\rho_1|>0
 \ \
 \mbox{whenever}
 \ \
 0<|\rho_1|<\sup_{n\geq 2}|\rho_n|.
 \]

 \noindent
 {\small\bf Remark.}
 %In the case where $f_1$ and $f_2$
 %are standard uniform densities,
 (a) It is obvious that the construction
 (\ref{eq.density}) can
 be adapted to the discrete (lattice) case where
 $(X,Y)\in\{1,\ldots, N\}^2$,
 covering the characterizations (for finite populations)
 treated by L\'{o}pez-Bl\'{a}zquez and Casta\~{n}o-Mart\'{i}nez
 (2006) and Casta\~{n}o-Mart\'{i}nez, L\'{o}pez-Bl\'{a}zquez
 and Salamanca-Mi\~{n}o (2007).
 \smallskip

 \noindent
 (b) Distributions with densities of the form
 (\ref{eq.density})
 %entails a generalization of the so called
 are known as Lancaster distributions;
 see, e.g.,
 %Lancaster (1957, 1969),
 Koudou (1998) or Diaconis and
 Griffiths (2012).
 They can be viewed as extensions
 of the Sarmanov-type distribution ($\rho_n=0$ for $n\geq 2$)
 which, assuming standard uniform marginals, generalizes
 the so called Farlie-Gumbel-Morngestern family.
 \vspace{1.5em}

 \noindent
 {\small\bf Acknowledgement.} I would like to thank
 H.N.\ Nagaraja and H.M.\ Barakat for bringing to my
 attention the papers by Rohatgi and Sz\'{e}kely (1992)
 and Sarmanov (1958a, 1958b).

 \vspace{2em}

 %\newpage
 \begin{center}
 {{\bf REFERENCES} }
 \end{center}

 %{ \vspace{1em} }
 {\small

 \begin{description}

 \item Afendras, G.; Papadatos, N.; Papathanasiou, V.\
 (2011). An extended Stein-type covariance identity for
 the Pearson family, with applications to lower variance
 bounds.
 {\it Bernoulli}, {\bf 17}, 507--529.
 \vspace{-.7ex}

 \item Ahsanullah, M.\ (2004). {\it Record Values -- Theory
 and Applications}. Univ.\ Press Amer.\ Inc., New York.
 \vspace{-.7ex}

 \item Arnold, B.C.; Balakrishnan, N.; Nagaraja, H.N.\
 (1998). {\it Records}. Wiley, New York.
 \vspace{-.7ex}

 \item Barakat, H.M.\ (2012). The maximal correlation for the
 generalized order statistics and dual generalized order
 statistics.
 {\it Arab.\ J.\ Math.}, {\bf 1}, 149--158.
 \vspace{-.7ex}

 %\item Balakrishnan, N.; Charalambides, C.; Papadatos, N.\
 %(2003). Bounds on expectation of order statistics from a finite
 %population,
 %{\it J.\ Statist.\ Plann.\ Inference}, {\bf 113}, 569--588.
 %\vspace{-.7ex}

 \item
 Berg, C.; Christensen, J.P.R.\ (1981).
 Density questions in the classical theory of moments.
 {\it Ann.\ Inst.\ Fourier} (Grenoble), {\bf 31},
 99--114.
 %\MR{0638619}
 \vspace{-.7ex}

 %\item Breiman, L.; Friedman, J.H.\ (1985). Estimating
 %optimal transformations for multiple regression and
 %correlation. {\it J.\ Amer.\ Statist.\ Assoc.},
 %{\bf 80}, 590--598.
 %\vspace{-.7ex}

 %\item Bryc, W.; Dembo, A.; Kagan, A.\
 %(2005). On the maximum correlation coefficient.
 %{\it Theory Probab.\ Appl.}, {\bf 49}, 132--138.
 %\vspace{-.7ex}

 \item Casta\~{n}o-Mart\'{i}nez, A.; L\'{o}pez-Bl\'{a}zquez, F.;
 Salamanca-Mi\~{n}o, B.\ (2007).
 Maximal correlation between order statistics.
 In: {\it Recent Developments in Ordered Random Variables},
 M.\ Ahsanullah and M.\ Raqab (eds.), Nova Science
 Publishers, 55--68.
 \vspace{-.7ex}

 \item Cs\'{a}ki, P.; Fischer, J.\ (1963). On the general
 notion of maximal correlation.
 {\it Magyar Tudom\'{a}nyos
 Akad.\ Mat. Kutat\'{o} Int\'{e}zetenk K\"{o}zlem\'{e}nyei}
 (Publ.\ Math.\ Inst.\ Hungar.\ Acad.\ Sci.),
 {\bf 8}, 27--51.
 \vspace{-.7ex}

 \item David, H.A.; Nagaraja, H.N.\ (2003). {\it Order
 Statistics}. Wiley, New York.
 \vspace{-.7ex}

 \item Diaconis, P.; Griffiths, R.\ (2012). Exchangeable
 pairs of Bernoulli random variables, Krawtchouck
 polynomials, and Ehrenfest urns. {\it
 Austral.\ New Zeal.\ J.\ Statist.}, (to appear).
 DOI: 10.1111/j.1467-842X.2012.00654.x
 \vspace{-.7ex}

 \item Dembo, A.; Kagan, A.; Shepp., L.A.\
 (2001). Remarks on the maximum correlation coefficient.
 {\it Bernoulli}, {\bf 7}, 343--350.
 \vspace{-.7ex}

 \item Eaton, M.\ (1986). A characterization of spherical
 distributions. {\it J.\ Multivariate Anal.}, {\bf 20},
 272--276.
 \vspace{-.7ex}

 \item Gebelein, H.\ (1941). Das Statistische Problem der
 Korrelation als Variation
 und Eigenwertproblem und sein Zusammenhang mit
 der Ausgleichrechnung.
 {\it Angew.\ Math.\ Mech.}, {\bf 21}, 364--379.
 \vspace{-.7ex}

 \item Hirschfeld, H.O.\ (1935). A connection between
 correlation and contingency. {\it Proc.\ Cambridge Phil.\
 Soc.}, {\bf 31}, 520--524.
 \vspace{-.7ex}

 %\item Huang, J.S.; Kotz, S.\ (1984). Correlation structure
 %of the iterated Farlie-Gumbel-Morgenstern distributions.
 %{\it Biometrika}, {\bf 71}(3), 633--636.
 %\vspace{-.7ex}

 %\item Koudou, A.E.\ (1996). Probabilit\'{e}s de Lancaster.
 %{\it Exp.\ Math.}, {\bf 14}, 247--275.
 %\vspace{-.7ex}

 \item Koudou, A.E.\ (1998). Lancaster bivariate probability distributions
 with Poisson, negative binomial and gamma margins.
 {\it Test}, {\bf 7}, 95--110.
 \vspace{-.7ex}

 %\item Lancaster, H.O.\ (1957). Some properties of the
 %bivariate normal distribution considered in the form of
 %a contingency table. {\it Biometrika}, {\bf 44}, 289--292.
 %\vspace{-.7ex}

 %\item Lancaster, H.O.\ (1969).
 %{\it The chi-squared distribution}. Wiley, New York.
 %\vspace{-.7ex}

 %\item Liu, J.S.; Wong, W.H.; Kong, A.\ (1994). Covariance
 %structure of the Gibbs sampler with applications to the
 %comparisons of estimators and augmentation schemes.
 %{\it Biometrika}, {\bf 81}, 27--40.
 %\vspace{-.7ex}

 \item L\'{o}pez-Bl\'{a}zquez, F.; Casta\~{n}o-Mart\'{i}nez, A.\
 (2006). Upper and lower bounds for the correlation ratio of order
 statistics from a sample without replacement.
 {\it J.\
 Statist.\ Plann.\ Inference}, {\bf 136}, 43--52.
 \vspace{-.7ex}

 \item Nevzorov, V.B.\ (1992). A characterization of
 exponential distributions
 by correlations between records.
 {\it  Math.\ Meth.\ Statist.}, {\bf 1}, 49--54.
 \vspace{-.7ex}

 %\item
 %Nomokonov, M.K.\ (1950). On the simpleness of the second
 %characteristic value
 %of correlation integral equations.
 %{\it Dokl.\ Akad.\ Nauk SSSR},
 %{\bf 72}, 1021--1024 (in
 %\vspace{-.7ex}
 %Russian).

 \item Papadatos, N.; Xifara, T.\ (2012). A simple method for
 obtaining the
 maximal correlation coefficient and related
 characterizations, {\it submitted for publication}.
 \vspace{-.7ex}
 arXiv:1204.1632v2.

 \item
 R\'{e}nyi, A.\ (1959).
 On measures of dependence.
 {\it Acta Math.\ Acad.\ Sci.\ Hungar.},
 {\bf 10}, 441-–451.
 \vspace{-.7ex}

 \item
 Rohatgi, V.K.; Sz\'{e}kely, G.J.\ (1992).
 On the background of some correlation inequalities.
 {\it J.\ Stat.\ Comput.\ Simul.}, {\bf 40}, 260-–262.
 \vspace{-.7ex}

 %\item
 %Sarmanov, O.V.\ (1946). On monotonic solutions
 %of correlation integral equations.
 %{\it Dokl.\ Akad.\ Nauk SSSR},
 %{\bf 53}, 781--784 (in Russian).
 %\vspace{-.7ex}

 \item
 Sarmanov, O.V.\ (1958a). Maximum correlation coefficient
 (symmetric case). {\it Dokl.\ Akad.\ Nauk SSSR},
 {\bf 120}, 715--718 (in Russian). Also
 in {\it Selected Translations in Mathematical Statistics
 and Probability}, {\bf 4} (1963),
 {\it Amer.\ Math.\ Soc.}, pp.\ 271--275.
 \vspace{-.7ex}

 \item
 Sarmanov, O.V.\ (1958b). Maximum correlation
 coefficient (non-symmetric case). {\it
 Dokl.\ Akad.\ Nauk SSSR},
 {\bf 121}, 52--55 (in Russian). Also
 in {\it Selected Translations in Mathematical Statistics
 and Probability}, {\bf 2} (1962),
 {\it Amer.\ Math.\ Soc.}, pp.\ 207--210.
 \vspace{-.7ex}

 \item Sz\'{e}kely, G.J.; Gupta, A.K.\ (1998).
 On a paper of V.B.\ Nevzorov.
 {\it Math.\ Meth.\ Statist.}, {\bf 7},
 122.
 \vspace{-.7ex}

 \item Sz\'{e}kely, G.J.; M\'{o}ri, T.F.\ (1985).
 An extremal property of rectangular
 distributions.
 {\it Statist.\ Probab.\ Lett.}, {\bf 3},
 107--109.
 \vspace{-.7ex}

 \item Terrell, G.R.\ (1983). A characterization of rectangular
 distributions.
 {\it Ann.\ Probab.}, {\bf 11}, 823--826.
 \vspace{-.7ex}

 %\item Yu, Y.\ (2008). On the maximal correlation
 %coefficient.
 %{\it Statist.\ Probab.\ Lett.}, {\bf 78}, 1072--1075.
 %\vspace{-.7ex}

 \item Vershik, A.M.\ (1964).
 Some characteristic properties of Gaussian stochastic processes.
 {\it Theory Probab.\ Appl.}, {\bf 9},
 390--394.
 %353--356 (in Russian).

 \end{description}
 }

 \end{document}